\documentclass{elsart}
\usepackage[latin1]{inputenc}
\usepackage[english]{babel}
\usepackage{amssymb}
\usepackage{amsfonts}
\usepackage{amsmath}
\usepackage{natbib}
\usepackage{url}
\usepackage{fancyhdr}
\usepackage{graphicx}
\usepackage{calc}

\newcommand{\R}[1]{\ensuremath{\Rset^{#1}}}

\newcommand{\intd}{\mathrm{d}}
\newcommand{\classifier}{f}

\newcommand{\sign}[1]{\textrm{sign}\left(#1\right)}
\newcommand{\indic}[1]{\mathbb{I}_{\left\{#1\right\}}}

\newtheorem{theorem}{Theorem}
\theoremstyle{remark}
\newtheorem{remark}{Remark}

\begin{document}
\begin{frontmatter}
\title{Support Vector Machine For Functional Data Classification}
\author[INRIA]{Fabrice Rossi\corauthref{rossi}},
\ead{Fabrice.Rossi@inria.fr}
\author[GRIMM]{Nathalie Villa},
\ead{villa@univ-tlse2.fr}
\address[INRIA]{Projet AxIS, INRIA-Rocquencourt, Domaine de Voluceau, 
  Rocquencourt, B.P.~105, 78153 Le Chesnay Cedex, France}
\address[GRIMM]{ Equipe GRIMM - Université Toulouse Le Mirail,
  5 allées A. Machado, 31058 Toulouse cedex 1 - FRANCE}
\corauth[rossi]{Corresponding author:\\
Fabrice Rossi\\
Projet AxIS\\
INRIA Rocquencourt\\
Domaine de Voluceau, Rocquencourt, B.P. 105\\
78153 LE CHESNAY CEDEX -- FRANCE\\
Tel: (33) 1 39 63 54 45\\
Fax: (33) 1 39 63 58 92\\
}
\end{frontmatter}

\newpage

\begin{frontmatter}
\title{Support Vector Machine For Functional Data Classification}
\begin{abstract}
  In many applications, input data are sampled functions taking their values
  in infinite dimensional spaces rather than standard vectors. This fact has
  complex consequences on data analysis algorithms that motivate modifications
  of them. In fact most of the traditional data analysis tools for regression,
  classification and clustering have been adapted to functional inputs under
  the general name of Functional Data Analysis (FDA). In this paper, we
  investigate the use of Support Vector Machines (SVMs) for functional data
  analysis and we focus on the problem of curves discrimination. SVMs are
  large margin classifier tools based on implicit non linear mappings of the
  considered data into high dimensional spaces thanks to kernels. We show how
  to define simple kernels that take into account the functional nature of the
  data and lead to consistent classification. Experiments conducted on real
  world data emphasize the benefit of taking into account some functional
  aspects of the problems.
\end{abstract}
\begin{keyword}
Functional Data Analysis\sep Support Vector Machine\sep Classification\sep Consistency
\end{keyword}
\end{frontmatter}

\newpage

\section{Introduction}\label{intro}
In many real world applications, data should be considered as discretized
functions rather than as standard vectors. In these applications, each
observation corresponds to a mapping between some conditions (that might be
implicit) and the observed response. A well studied example of those
functional data is given by spectrometric data (see section
\ref{app_spectro}): each spectrum is a function that maps the wavelengths of
the illuminating light to the corresponding absorbances (the responses) of the
studied sample. Other natural examples can be found in voice recognition area
(see sections \ref{app_biau} and \ref{app_hastie}) or in meteorological
problems, and more generally, in multiple time series analysis where each
observation is a complete time series.

The direct use of classical models for this type of data faces several
difficulties: as the inputs are discretized functions, they are generally
represented by high dimensional vectors whose coordinates are highly
correlated. As a consequence, classical methods lead to ill-posed problems,
both on a theoretical point of view (when working in functional spaces that
have infinite dimension) and on a practical one (when working with the
discretized functions). The goal of Functional Data Analysis (FDA) is to use,
in data analysis algorithms, the underlying functional nature of the data:
many data analysis methods have been adapted to functions (see
\cite{RamseySilverman97} for a comprehensive introduction to functional data
analysis and a review of linear methods). While the original papers on FDA
focused on linear methods such as Principal Component Analysis
\cite{Deville74,DauxoisPousse76,DauxoisPousseRomain82,BesseRamsay1986} and the
linear model \cite{RamseyDalzell1991,FrankFriedman1993,HastieMallows1993}, non
linear models have been studied extensively in the recent years. This is the
case, for instance, of most neural network models
\cite{FerreVilla04SIRNN,RossiConanGuez05NeuralNetworks,RossiConanGuezElGolliESANN2004SOMFunc,RossiEtAl05Neurocomputing}.

In the present paper, we adapt Support Vector Machines (SVMs, see e.g.
\cite{Vapnik1995,ChristianiniShaweTaylor2000SVMIntroduction}) to functional
data classification (the paper extends results from
\cite{RossiVillaASMDA2005SVM,VillaRossiESANN2005SVM}). We show in particular
both the practical and theoretical advantages of using functional kernels,
which are kernels that take into account the functional nature of the data. On
a practical point of view, those kernels allow to take advantage of the expert
knowledge on the data. On the theoretical point of view, a specific type of
functional kernels allows the construction of a consistent training procedure
for functional SVMs.

The paper is organized as follow: section \ref{FDA} presents the functional
data classification and why it generally leads to ill-posed problems. Section
\ref{SVM fonct} provides a short introduction to SVMs and explains why their
generalization to FDA can lead to particular problems. Section
\ref{sectionFDAKernels} describes several functional kernels and explains how
they can be practically computed while section \ref{sectionConsistency}
presents a consistency result for some of them. Finally, section
\ref{applications} illustrates the various approaches presented in the paper
on real data sets.

\section{Functional Data Analysis}
\label{FDA}
\subsection{Functional Data}\label{subsectionFDA}
To simplify the presentation, this article focuses on functional data for
which each observation is described by one function from \R{} to \R{}.
Extension to the case of several real valued functions is straightforward.
More formally, if $\mu$ denotes a known finite positive Borel measure on $\R{}$, an
observation is an element of $L^2(\mu)$, the Hilbert space of
$\mu$-square-integrable real valued functions defined on $\R{}$. In some
situations, additional regularity assumptions (e.g., existence of derivatives)
will be needed.

However, almost all the developments of this paper are not specific to
functions and use only the Hilbert space structure of $L^2(\mu)$. We will
therefore denote $\mathcal{X}$ an arbitrary Hilbert space and $\langle
.,.\rangle$ the corresponding inner product.  Additional assumptions on
$\mathcal{X}$ will be given on a case by case basis.  As stated above, the
most common situation will of course be $\mathcal{X}=L^2(\mu)$ with $\langle
u,v\rangle=\int uv\intd\mu$.

\subsection{Data analysis methods for Hilbert spaces}%
\label{subsectionmethod} 
It should be first noted that many data analysis algorithms can be written so
as to apply, at least on a theoretical point of view, to arbitrary Hilbert
spaces. This is obviously the case, for instance, for distance-based
algorithms such as the $k$-nearest neighbor method. Indeed, this algorithm
uses only the fact that distances between observations can be calculated.
Obviously, it can be applied to Hilbert spaces using the distance induced by
the inner product. This is also the case of methods directly based on inner
products such as multi-layer perceptrons (see
\cite{Sandberg1996,SandbergXu1996,Stinchcombe99} for a presentation of
multi-layer perceptrons with almost arbitrary input spaces, including Hilbert
spaces).

However, functional spaces have infinite dimension and a basic transposition
of standard algorithms introduces both theoretical and practical difficulties.
In fact, some simple problems in \R{d} become ill-posed in $\mathcal{X}$ when
the space has infinite dimension, even on a theoretical point of view.

Let us consider for instance the linear regression model in which a real
valued target variable $Y$ is modeled by $E(Y|X)=H(X)$ where $H$ is a linear
continuous operator defined on the input space. When $X$ has values in \R{d}
(i.e., $\mathcal{X}=\R{d}$), $H$ can be easily estimated by the least
square method that leads to the inversion of the covariance matrix of $X$. In
practice, problems might appear when $d$ is not small compared to $N$, the
number of available examples, and regularization techniques should be used
(e.g., ridge regression \cite{HoerlRidge1970}). When $X$ has values in a
Hilbert space, the problem is ill-posed because the covariance of $X$
is a Hilbert-Schmidt operator and thus has no continuous inverse; direct approximation of the inverse of this operator is then problematic as it does not provide a consistant estimate (see \cite{CardotFerratySarda1999}).

To overcome the infinite dimensional problem, most of FDA methods so far have
been constructed thanks to two general principles: \emph{filtering} and
\emph{regularization}. In the filtering approach, the idea is to use
representation methods that allow to work in finite dimension (see for
instance \cite{CardotFerratySarda1999} for the functional linear model and
\cite{BiauEtAl2005FunClassif} for a functional $k$-nearest neighbor method).
In the regularization approach, the complexity of the solution is constrained
thanks to smoothness constraints. For instance, building a linear model in a
Hilbert space consists in finding a function $h\in L^2(\mu)$ such that
$E(Y|X)=\langle h,X\rangle$. In the regularization approach, $h$ is chosen
among smooth candidates (for instance twice derivable functions with minimal
curvature), see e.g.
\cite{HastieMallows1993,MarxEilers1996,CardotFerratySarda2002}. Other examples
of the regularization approach include smooth Principal Component Analysis
\cite{PezzulliSilverman1993} and penalized Canonical Component Analysis
\cite{LeurgansMoyeedSilverman1993}. A comparison of filtering and
regularization approaches for a semi-parametric model used in curve
discrimination can be found in \cite{ferre_villa_RSA2005}.

Using both approaches, a lot of data analysis algorithms have been
successfully adapted to functional data. Our goal in the present paper is to
study the case of Support Vector Machines (SVM), mainly thanks to a filtering
approach.

\section{Support Vector Machines for FDA}
\label{SVM fonct}
\subsection{Support Vector Machines}\label{subsectionSVM}
We give, in this section, a very brief presentation of Support Vector Machines
(SVMs) that is needed for the definition of their functional versions. We
refer the reader to e.g. \cite{ChristianiniShaweTaylor2000SVMIntroduction} for
a more comprehensive presentation. As stated in section~\ref{subsectionFDA},
$\mathcal{X}$ denotes an arbitrary Hilbert space. Our presentation of SVM
departs from the standard introduction because it assumes that the
observations belong to $\mathcal{X}$ rather than to a \R{d}. This will make
clear that the definition of SVM on arbitrary Hilbert spaces is not the
difficult part in the construction of functional SVM. We will discuss problems
related to the functional nature of the data in section
\ref{subsectionSVMFDA}.

Our goal is to classify data into two predefined classes. We assume given a
learning set, i.e. $N$ examples $(x_1,y_1),\ldots,(x_N,y_N)$ which are i.i.d.
realizations of the random variable pair $(X,Y)$ where $X$ has values in
$\mathcal{X}$ and $Y$ in $\{-1,1\}$, i.e. $Y$ is the class label for $X$ which
is the observation. 

\subsubsection{Hard margin SVM}
The principle of SVM is to perform an affine discrimination of the
observations with maximal margin, that is to find an element $w\in
\mathcal{X}$ with a minimum norm and a real value $b$, such that $y_i(\langle
w,x_i\rangle + b)\geq 1$ for all $i$. To do so, we have to solve the following
quadratic programming problem:
\[
(P_0)\ \min_{w,b} \langle w,w\rangle ,\textrm{ subject to } y_i(\langle w,x_i\rangle +b)\geq 1,\ 1\leq i\leq N.
\]
The classification rule associated to $(w,b)$ is simply
$\classifier(x)=\mathrm{sign}(\langle w,x\rangle + b)$. In this situation (called hard
margin SVM), we request the rule to have zero error on the learning set.

\subsubsection{Soft margin SVM}\label{subsubsectionSoftMarginSVM}
In practice, the solution provided by problem $(P_0)$ is not very
satisfactory. Firstly, perfectly linearly separable problems are quite rare,
partly because non linear problems are frequent, but also because noise can
turn a linearly separable problem into a non separable one. Secondly, 
choosing a classifier with maximal margin does not prevent overfitting,
especially in very high dimensional spaces (see
e.g. \cite{HastieEtAl2004EntireSVM} for a discussion about this point). 

A first step to solve this problem is to allow some classification errors on
the learning set. This is done by replacing $(P_0)$ by its soft margin
version, i.e., by the problem:
\[
(P_C) \begin{array}[t]{l}
\min_{w,b,\xi} \langle w,w\rangle +C\sum_{i=1}^N\xi_i,\\
\textrm{subject to } \begin{array}[t]{l}
y_i(\langle w,x_i\rangle +b)\geq 1-\xi_i,\ 1\leq i\leq N,\\
\xi_i\geq 0,\ 1\leq i\leq N.
\end{array}
\end{array}
\]
Classification errors are allowed thanks to the slack variables $\xi_i$.  The
$C$ parameter acts as an inverse regularization parameter. When $C$ is small,
the cost of violating the hard margin constraints, i.e., the cost of having
some $\xi_i>0$ is small and therefore the constraint on $w$ dominates. On the
contrary, when $C$ is large, classification errors dominate and $(P_C)$ gets
closer to $(P_0)$. 

\subsubsection{Non linear SVM}
As noted in the previous section, some classification problems don't have a
satisfactory linear solution but have a non linear one. Non linear SVMs are
obtained by transforming the original data. Assume given an
Hilbert space $\mathcal{H}$ (and denote $\langle
.,.\rangle_{\mathcal{H}}$ the corresponding inner product) and a function
$\phi$ from $\mathcal{X}$ to $\mathcal{H}$ (this function is called a
\emph{feature map}). A linear SVM in $\mathcal{H}$ can be constructed on the
data set $(\phi(x_1),y_1),\ldots,(\phi(x_N),y_N)$. If $\phi$ is a non linear
mapping, the classification rule $\classifier(x)=\mathrm{sign}(\langle
w,\phi(x)\rangle_{\mathcal{H}} + b)$ is also non linear.

In order to obtain the linear SVM in $\mathcal{H}$ one has to solve the
following optimization problem:
\[
(P_{C,\mathcal{H}}) \begin{array}[t]{l}
\min_{w,b,\xi} \langle w,w\rangle_{\mathcal{H}} +C\sum_{i=1}^N\xi_i,\\
\textrm{subject to } \begin{array}[t]{l}
y_i(\langle w,\phi(x_i)\rangle_{\mathcal{H}} +b)\geq 1-\xi_i,\ 1\leq i\leq N,\\
\xi_i\geq 0,\ 1\leq i\leq N.
\end{array}
\end{array}
\]
It should be noted that this feature mapping allows to define SVM on almost
arbitrary input spaces. 

\subsubsection{Dual formulation and Kernels}
Solving problems $(P_C)$ or $(P_{C,\mathcal{H}})$ might seem very difficult at
first, because $\mathcal{X}$ and $\mathcal{H}$ are arbitrary Hilbert spaces
and can therefore have very high or even infinite dimension (when
$\mathcal{X}$ is a functional space for instance). However, each problem has a
dual formulation. More precisely, $(P_{C})$ is equivalent to the
following optimization problem (see \cite{Lin2001SVM}):
\[
(D_{C}) \begin{array}[t]{l}
\max_{\alpha} \sum_{i=1}^N\alpha_i-\sum_{i=1}^N\sum_{j=1}^N\alpha_i \alpha_j y_i y_j \langle x_i,x_j\rangle ,\\
\textrm{subject to } \begin{array}[t]{l}
\sum_{i=1}^N\alpha_iy_i=0,\\
0\leq\alpha_i\leq C,\ 1\leq i\leq N.
\end{array}
\end{array}
\]
This result applies to the original problem in which data are not mapped into
$\mathcal{H}$, but also to the mapped data, i.e., $(P_{C,\mathcal{H}})$ is
equivalent to a problem $(D_{C,\mathcal{H}})$ in which the $x_i$ are replaced
by $\phi(x_i)$ and in which the inner product of $\mathcal{H}$ is used. This
leads to: 
\[
(D_{C,\mathcal{H}}) \begin{array}[t]{l}
\max_{\alpha} \sum_{i=1}^N\alpha_i-\sum_{i=1}^N\sum_{j=1}^N\alpha_i \alpha_j y_i y_j \langle \phi(x_i),\phi(x_j)\rangle_{\mathcal{H}} ,\\
\textrm{subject to } \begin{array}[t]{l}
\sum_{i=1}^N\alpha_iy_i=0,\\
0\leq\alpha_i\leq C,\ 1\leq i\leq N.
\end{array}
\end{array}
\]
Solving $(D_{C,\mathcal{H}})$ rather than $(P_{C,\mathcal{H}})$ has two
advantages. The first positive aspect is that $(D_{C,\mathcal{H}})$ is an
optimization problem in $\R{N}$ rather than in $\mathcal{H}$ which can have
infinite dimension (the same is true for $\mathcal{X}$). 

The second important point is linked to the fact that the optimal
classification rule can be written
$\classifier(x)=\mathrm{sign}(\sum_{i=1}^N\alpha_iy_i\langle
\phi(x_i),\phi(x)\rangle_{\mathcal{H}} + b)$. This means that both the
optimization problem and the classification rule do not make direct use of the
transformed data, i.e. of the $\phi(x_i)$. All the calculations are done
through the inner product in $\mathcal{H}$, more precisely through the
values $\langle \phi(x_i),\phi(x_j)\rangle_{\mathcal{H}}$. Therefore, rather
than choosing directly $\mathcal{H}$ and $\phi$, one can provide a so called
\emph{Kernel function} $K$ such that $K(x_i,x_j)=\langle
\phi(x_i),\phi(x_j)\rangle_{\mathcal{H}}$ for a given pair
$(\mathcal{H},\phi)$. 

In order that $K$ corresponds to an actual inner product in a Hilbert space,
it has to fulfill some conditions. $K$ has to be symmetric and positive
definite, that is, for every $N$, $x_1,\ldots,x_N$ in $\mathcal{X}$ and
$\alpha_1,\ldots,\alpha_N$ in \R{},
$\sum_{i=1}^N\sum_{j=1}^N\alpha_i\alpha_jK(x_i,x_j)\geq 0$. If $K$ satisfies
those conditions, according to Moore-Aronszajn theorem \cite{Aronszajn1950},
there exists a Hilbert space $\mathcal{H}$ and feature map $\phi$ such that
$K(x_i,x_j)=\langle \phi(x_i),\phi(x_j)\rangle_{\mathcal{H}}$. 

\subsection{The case of functional data}\label{subsectionSVMFDA}
The short introduction to SVM proposed in the previous section has clearly
shown that defining linear SVM for data in a functional space is as easy as
for data in \R{d}, because we only assumed that the input space was a Hilbert
space. By the dual formulation of the optimization problem $(P_C)$, a
software implementation of linear SVM on functional data is even possible, by
relying on numerical quadrature methods to calculate the requested integrals
(inner product in $L^2(\mu)$, cf section \ref{sectionInpractice}). 

However, the functional nature of the data has some effects. It should be
first noted that in infinite dimensional Hilbert spaces, the hard margin
problem $(P_0)$ has always a solution when the input data are in general
positions, i.e., when $N$ observations span a $N$ dimensional subspace of
$\mathcal{X}$. A very naive solution would therefore consists in avoiding soft
margins and non linear kernels. This would not give very interesting results
in practice because of the lack of regularization (see
\cite{HastieEtAl2004EntireSVM} for some examples in very high dimension
spaces, as well as section \ref{app_biau}).

Moreover, the linear SVM with soft margin can also lead to bad
performances. It is indeed well known (see
e.g. \cite{HastieTibshiraniFriedman2001SL}) that problem $(P_C)$ is equivalent
to the following unconstrained optimization problem:
\[
(R_\lambda) \min_{w,b}\frac{1}{N}\sum_{i=1}^N\max\left(0,1-y_i(\langle w,x_i\rangle +b)\right)+\lambda\langle w,w\rangle,
\]
with $\lambda=\frac{1}{CN}$. This way of viewing $(P_C)$ emphasizes the
regularization aspect (see also
\cite{SmolaScholkpofAlgorithmica98,SmolaEtAlNN98,EvgeniouEtAlRegularization2000})
and links the SVM model to ridge regression \cite{HoerlRidge1970}. As shown in
\cite{HastieBujaTibshirani1995}, the penalization used in ridge regression
behaves poorly with functional data. Of course, the loss function used by SVM
(the \emph{hinge loss}, i.e., $h(u,v)=\max(0,1-uv)$) is different from the
quadratic loss used in ridge regression and therefore no conclusion can be
drawn from experiments reported in \cite{HastieBujaTibshirani1995}. However
they show that we might expect bad performances with the linear SVM applied
directly to functional data. We will see in sections \ref{app_biau} and
\ref{app_hastie} that the efficiency of the ridge regularization seems to be
linked with the actual dimension of the data: it does not behave very well
when the number of discretization points is very big and thus leads to
approximate the ridge penalty by a dot product in a very high dimensional
space (see also section \ref{sectionInpractice}).

It is therefore interesting to consider non linear SVM for functional data, by
introducing adapted kernels. As pointed out in
e.g. \cite{EvgeniouEtAlRegularization2000}, $(P_{C,\mathcal{H}})$ is equivalent
to 
\[
(R_{\lambda,\mathcal{H}}) \min_{f\in\mathcal{H}}\frac{1}{N}\sum_{i=1}^N\max\left(0,1-y_if(x_i))\right)+\lambda\langle f,f\rangle_{\mathcal{H}}.
\]
Using a kernel corresponds therefore both to replace a linear classifier by a
non linear one, but also to replace the ridge penalization by a penalization
induced by the kernel which might be more adapted to the problem (see
\cite{SmolaEtAlNN98} for links between regularization operators and kernels). The applications presented in \ref{applications} illustrate this fact.

\section{Kernels for FDA}\label{sectionFDAKernels}
\subsection{Classical kernels}
Many standard kernels for \R{d} data are based on the Hilbert structure of
\R{d} and can therefore be applied to any Hilbert space. This is the case for
instance of the Gaussian kernel (based on the norm in $\mathcal{X}$: $K(u,v)=e^{-\sigma \| u-v\|^2}$) and of the
polynomial kernels (based on the inner product in $\mathcal{X}$: $K(u,v)=(1+\langle u,v\rangle)^D$). Obviously, the only
practical difficulty consists  in implementing the calculations needed
in $\mathcal{X}$ so as to evaluate the chosen kernel (the problem also appears
for the plain linear ``kernel'', i.e. when no feature mapping is
done). Section \ref{sectionInpractice} discusses this point.

\subsection{Using the functional nature of the data}
While the functional version of the standard kernels can provide an
interesting library of kernels, they do not take advantage of the functional
nature of the data (they use only the Hilbert structure of
$L^2(\mu)$). Kernels that use the fact that we are dealing with functions are
nevertheless quite easy to define.

A standard method consists in introducing kernels that are made by a
composition of a simple feature map with a standard kernel. More formally, we
use a transformation operator $P$ from $\mathcal{X}$ to another space
$\mathcal{D}$ on which a kernel $K$ is defined. The actual kernel $Q$ on
$\mathcal{X}$ is defined as $Q(u,v)=K(P(u),P(v))$ (if $K$ is a kernel, then so
is $Q$). 

\subsubsection{Functional transformations}\label{sectionTransformation}
In some application domains, such as chemometrics, it is well known that the
shape of a spectrum (which is a function) is sometimes more important than its
actual mean value. Several transformations can be proposed to deal with this
kind of data. For instance, if $\mu$ is a finite measure (i.e.,
$\mu(\R{})<\infty$), a centering transformation can be defined as the
following mapping from $L^2(\mu)$ to itself:
\[
C(u)=u-\frac{1}{\mu(\R{})}\int u\intd \mu.
\]
A normalization mapping can also be defined:
\[
N(u)=\frac{1}{\|C(u)\|}C(u).
\]
If the functions are smooth enough, i.e., if we restrict ourselves to a
Sobolev space $W^{s,2}$, then some derivative transformations can be used: the
Sobolev space $W^{s,2}$, also denoted $H^s$, is the Hilbert space of functions
which have $L^2$ derivatives up to the order $s$ (in the sense of the
distribution theory). For instance, with $s\geq 2$, we can use the second
derivative that allows to focus on the curvature of the functions: this is
particularly useful in near infrared spectrometry (see e.g.,
\cite{RossiConanGuez05NeuralNetworks,RossiEtAl05Neurocomputing}, and section
\ref{app_spectro}). 

\subsubsection{Projections}\label{sectionProjection}
Another type of transformations can be used in order to define adapted
kernels. The idea is to reduce the dimensionality of the input space, that is
to apply the standard filtering approach of FDA. We assume given a
$d$-dimensional subspace $V_d$ of $\mathcal{X}$ and an orthonormal basis of this
space denoted $\{\Psi_j\}_{j=1,\ldots,d}$. We define the transformation $P_{V_d}$
as the orthogonal projection on $V_d$,
\[
P_{V_d}(x)=\sum_{j=1}^d \langle x,\Psi_j\rangle \Psi_j.
\]
$(V_d,\langle.,.\rangle_{\mathcal{X}})$ is isomorphic to
$(\R{d},\langle.,.\rangle_{\R{d}})$ and therefore one can use a standard \R{d}
SVM on the vector data $(\langle x,\Psi_1\rangle,\ldots,\langle
x,\Psi_d\rangle)$. This means that $K$ can be any kernel adapted to vector data. In the case where $K$ is the usual dot product of $\R{d}$, this kernel is known as the empirical kernel map (see \cite{vert_tsuda_scholkopf_KMICB2004} for further details in the field of protein analysis).

Obviously, this approach is not restricted to functional data, but the choice
of $V_d$ can be directed by expert knowledge on the considered functions and we
can then consider that it takes advantage of the functional nature of the
data. We outline here two possible solutions based on orthogonal basis and on
B-spline basis.

If $\mathcal{X}$ is separable, it has a Hilbert basis, i.e., a complete
orthonormal system $\{\Psi_j\}_{j\geq 1}$. Therefore one can define $V_d$ as
the space spanned by $\{\Psi_j\}_{j=1,\ldots,d}$. The choice of the basis can
be based on expert considerations. Good candidates include Fourier basis and
wavelet basis. If the signal is known to be non stationary, a wavelet based
representation might for instance give better results than a Fourier
representation. Once the basis is chosen, an optimal value for $d$ can be
derived from the data, as explained in section \ref{sectionConsistency}, in such a way
that the obtained SVM has some consistency properties. Moreover, this
projection approach gives good results in practice (see section
\ref{app_biau}).

Another solution is to choose a projection space that has interesting
practical properties, for instance a spline space with its associated B-spline
bases. Spline functions regularity can be chosen \emph{a priori} so as to
enforce expert knowledge on the functions. For instance, near infrared spectra
are smooth because of the physical properties of the light transmission (and
reflection). By using a spline representation of the spectra, we replace
original unconstrained observations by $C^k$ approximations ($k$ depends on
what kind of smoothness hypothesis can be done). This projection can also be
combined with a derivative transformation operation (as proposed in section
\ref{sectionTransformation}). 

\subsection{Functional data in practice}\label{sectionInpractice}
In practice, the functions $(x_i)_{1\leq i\leq N}$ are never perfectly
known. It is therefore difficult to implement exactly the functional kernels
described in this section. 

The best situation is the one in which $d$ discretization points have been
chosen in $\R{}$, $(t_k)_{1\leq k\leq d}$, and each function $x_i$ is described
by a vector of \R{d}, $\left(x_i(t_1),\ldots,x_i(t_d)\right)$. In this
situation, a simple solution consists in assuming that standard operations in
\R{d} (linear combinations, inner product and norm) are good approximations of
their counterparts in the considered functional space. When the sampling is
regular, this is equivalent to applying standard SVMs to the vector
representation of the functions (see section \ref{applications} for real world
examples of this situation). When the sampling is not regular, integrals
should be approximated thanks to a quadrature method that will take into
account the relative position of the sampling points.

In some application domains, especially medical ones (e.g.,
\cite{JamesHastie2001}), the situation is not as good. Each function is in
general badly sampled: the number and the location of discretization points
depend on the function and therefore a simple vector model is not anymore
possible. A possible solution in this context consists in constructing a
approximation of $x_i$ based on its observation values (thanks to e.g.,
B-splines) and then to work with the reconstructed functions (see
\cite{RamseySilverman97,RossiEtAl05Neurocomputing} for details). 

The function approximation tool used should be simple enough to allow easy
implementation of the requested operations. This is the case for instance for
B-splines that allow in addition derivative calculations and an easy
implementation of the kernels described in section
\ref{sectionTransformation}. It should be noted that spline approximation is
different from projection on a spline subspace. Indeed each sampled function
could be approximated on a different B-spline basis, whereas the projection
operator proposed in section \ref{sectionProjection} requests an unique
projection space and therefore the same B-spline basis for each input
function. In other words, the spline approximation is a convenient way of
representing functions (see section \ref{app_spectro} for an application to
real world data), whereas the spline projection corresponds to a data
reduction technique. Both aspects can be combined. 

\section{Consistency of functional SVM}\label{sectionConsistency}
\subsection{Introduction}
In this section we study one of the functional kernel described above and show
that it can be used to define a consistent classifier for functional data. We
introduce first some notations and definitions.

Our goal is to define a training procedure for functional SVM such that the
asymptotic generalization performances of the constructed model is optimal. We
define as usual the generalization error of a classifier $\classifier$
by the probability of misclassification:
\[
L\classifier = \mathbb{P}(\classifier(X)\neq Y).
\]
The minimal generalization error is the Bayes error achieved by the optimal
classifier $\classifier^*$  given by
\[
\classifier^*(x)=\left\{\begin{array}{cl} 1&\textrm{when }
    \mathbb{P}(Y=1\mid X=x)>1/2 \\ 
    -1 & \textrm{otherwise.}
\end{array}\right.
\]
We denote $L^*=L\classifier^*$ the optimal Bayes error. Of course, the closer
the error of a classifier is from $L^*$, the better its generalization ability
is. 

Suppose that we are given a learning sample of size $N$ defined as in section~\ref{subsectionSVM}. A learning procedure is an algorithm which allows the
construction, from this learning sample, of a classification rule
$\classifier_N$ chosen in a set of admissible classifiers. This
algorithm is said to be consistent if 
\[
L\classifier_N\xrightarrow{N\rightarrow+\infty}L^*. 
\]
It should be noted that when the data belong to \R{d}, SVMs don't always
provide consistent classifiers. Some sufficient conditions have been given in
\cite{SteinwartJC2002}: the input data must belong to a compact subset of
\R{d}, the regularization parameter ($C$ in $(P_{C,\mathcal{H}})$) has to be
chosen in specific way (in relation to $N$ and to the type of kernel used) and the kernel must be
\emph{universal} \cite{SteinwartJMLR2001}.  If $\phi$ is the feature map
associated to a kernel $K$, the kernel is universal if the set of all the
functions of the form $x\mapsto \langle w,\phi(x)\rangle$ for
$w\in\mathcal{H}$ is dense in the set of all continuous functions defined on
the considered compact subset. In particular, the Gaussian kernel with any
$\sigma>0$ is universal for all compact subsets of $\mathbb{R}^d$ (see \cite{SteinwartJC2002} for futher details and the proof of Theorem \ref{th consist} for the precise statement on $C$).

\subsection{A learning algorithm for functional SVM}\label{subsectionLearningSVM}
The general methodology proposed in \cite{BiauEtAl2005FunClassif} allows to
turn (with some adaptations) a consistent algorithm for data in \R{d} into a
consistent algorithm for data in $\mathcal{X}$, a separable Hilbert space. We
describe in this section the adapted algorithm based on SVM.

The methodology proposed in \cite{BiauEtAl2005FunClassif} is based on
projection operators described in section \ref{sectionProjection}, more
precisely on the usage of a Hilbert basis of $\mathcal{X}$. In order to build
a SVM classifier based on $N$ examples, one need to choose from the data
several parameters (in addition to the weights $\{\alpha_i\}_{1\leq i\leq N}$
and $b$ in problem $(D_{C,\mathcal{H}})$):
\begin{enumerate}
\item the projection size parameter $d$, i.e., the dimension of the subset
  $V_d$ on which the functions are projected before being submitted to the SVM
  (recall that $V_d$ is the space spanned by $\{\Psi_j\}_{j=1,\ldots,d}$);
\item $C$, the regularization parameter;
\item the fully specified kernel $K$, that is the type of the universal kernel
  (Gaussian, exponential, etc.)  but also the parameter of this kernel such as
  $\sigma$ for the Gaussian kernel $K(u,v)=e^{-\sigma^2\|u-v\|^2}$. 
\end{enumerate}
Let us denote $\mathcal{A}$ the set of lists of parameters to explore (see
section \ref{subsectionConsisency} for practical examples). Following
\cite{BiauEtAl2005FunClassif} we use a validation approach to choose the best
list of parameters $a\in\mathcal{A}$ and in fact the best classifier on the validation set.

The data are split into two sets: a training set $\{(x_i,y_i),
i=1,\ldots,l_N\}$ and a validation set $\{(x_i,y_i), i=l_N+1,\ldots,N\}$. For
each fixed list $a$ of parameters, the training set $\{(x_i,y_i),
i=1,\ldots,l_N\}$ is used to calculate the SVM classification rule
$\classifier_a(x)=\sign{\sum_{i=1}^{l_N} \alpha_i^* y_i
  K(P_{V_d}(x),P_{V_d}(x_i))+b^*}$ where $(\{\alpha_i^*\}_{1\leq i\leq
  l_N},b^*)$ is the solution of $(D_{C,\mathcal{H}})$ applied to the projected
data $\{P_{V_d}(x_i), i=1,\ldots,l_N\}$ (please note that everything should be
indexed by $a$, for instance one should write $K_a$ rather than $K$).

The validation set is used to select the optimal value of $a$ in
$\mathcal{A}$, $a^*$, according to estimation of the generalization error
based on a penalized empirical error, that is, we define
\[
a^*=\arg\min_{a\in\mathcal{A}}{\widehat{L}\classifier_a+\frac{\lambda_a}{\sqrt{N-l_N}}},
\]
where
\[
\widehat{L}\classifier_a=\frac{1}{N-l_N}\sum_{n=l_N+1}^N \indic{\classifier_a(x_n)\neq y_n},
\]
and $\lambda_a$ is a penalty term used to avoid the selection of the most
complex models (i.e., the one with the highest $d$ in general). The classifier
$\classifier_N$ is then chosen as $\classifier_N=\classifier_{a^*}$. 

\subsection{Consistency}\label{subsectionConsisency}
Under some conditions on $\mathcal{A}$, the algorithm proposed in the
previous section is consistent. We assume given a fixed Hilbert basis of the
separable Hilbert space $\mathcal{X}$, $\{\Psi_j\}_{j\geq 1}$. When the
dimension of the projection space $V_d$ is chosen, a fully specified kernel
$K$ has to be chosen in a finite set of kernels, $\mathcal{J}_d$. The regularization parameter $C$
can be chosen in a bounded interval of the form
$[0,\mathcal{C}_d]$, for instance thanks to the algorithm
proposed in \cite{HastieEtAl2004EntireSVM} that allows to calculate the
validation performances for all values of $C$ in a finite time. Therefore, the
set $\mathcal{A}$ can be written $\bigcup_{d\geq 1}\{d\}\times \mathcal{J}_d
\times [0,\mathcal{C}_d]$. An element of $\mathcal{A}$ is a triple $a=(d,K,C)$
that specifies the projection operator $P_{V_d}$, the kernel $K$ (including
all its parameters) and the regularization constant $C$.

Let us first define, for all $\epsilon > 0$,
$\mathcal{N}(\mathcal{H},\epsilon)$ the covering number of the Hilbert space
$\mathcal{H}$ which is the minimum number of balls with radius $\epsilon$ that
are needed to cover the whole space $\mathcal{H}$ (see e.g., chapter 28 of
\cite{DevroyeEtAl1996Pattern}). Note that in SVM, as $\mathcal{H}$ is induced
by a kernel $K$, this number is closely related to the kernel (in particular
because the norm used to defined the balls is induced by the inner product of
$\mathcal{H}$, that is by $K$ itself); in this case, we will then denote the
covering number $\mathcal{N}(K,\epsilon)$. For example, Gaussian kernels are
known to induce feature spaces with covering number of the form
$\mathcal{O}(\epsilon^{-d})$ where $d$ is the dimension of the input space
(see \cite{SteinwartJC2002}).

Then we have:
\begin{theorem}\label{th consist}
  We assume that $X$ takes its values in a bounded subspace of
  the separable Hilbert space $\mathcal{X}$. We suppose that,
  \begin{align*}
  \forall d \geq1,\qquad &
  \mathcal{J}_d\textrm{ is a finite set,}&
 \\
 & \exists K_d \in \mathcal{J}_d\textrm{ such that: }\begin{array}[t]{l}
 K_d\textrm{ is universal,}\\ \exists \nu_d>0:\ \mathcal{N}(K_d,\epsilon)=\mathcal{O}(\epsilon^{-\nu_d}),\end{array}&\\
 & \mathcal{C}_d> 1,&
\end{align*}
and that
\[
\sum_{d\geq1} |\mathcal{J}_d| e^{-2\lambda_d^2}<+\infty,
\]
and finally that
  \begin{align*}
\lim_{N\rightarrow+\infty}l_N=+\infty&&
\lim_{N\rightarrow +\infty} N-l_N=+\infty\\
\lim_{N\rightarrow +\infty} \frac{l_N \log(N-l_N)}{N-l_N}=0.&&
\end{align*}
Then, the functional SVM $\classifier_N=\classifier_{a^*}$ chosen as described
in section \ref{subsectionLearningSVM} (where $a^*$ is optimal in
$\mathcal{A}=\bigcup_{d\geq 1}\{d\}\times \mathcal{J}_d
\times [0,\mathcal{C}_d]$) is consistent that is:
\[
L\classifier_N\xrightarrow{N\rightarrow+\infty}L^*.
\]
\end{theorem}
The proof of this result is given in Appendix \ref{preuves}. It is close to
the proof given in \cite{BiauEtAl2005FunClassif} except that in
\cite{BiauEtAl2005FunClassif} the proof follows from an oracle inequality
given for a finite grid search model. The grid search is adapted to the
classifier used in the paper (a $k$-nearest neighbor method), but not to our
setting. Our result includes the search for a parameter $C$ which can belong
to an infinite and non countable set; this can be done by the use of the
shatter coefficient of a particular class of linear classifiers which provides
the behavior of the classification rule on a set of $N-l_N$ observations (see
\cite{DevroyeEtAl1996Pattern}).

As pointed out before, the Gaussian kernel satisfies the hypothesis of the
theorem. Therefore, if $\mathcal{I}_d$ contains a Gaussian kernel for all $d$,
then consistency of the whole procedure is guaranteed. Other non universal
kernels can of course be included in the search for the optimal model.

\begin{remark}
  Note that, in this theorem, the sets $\mathcal{J}_d$ and $[0,\mathcal{C}_d]$
  depend on $d$: this does not influence the consistency of the method. In
  fact, one could have chosen the same set for every $d$, and $\mathcal{J}_d$
  could also contain a single Gaussian kernel with any parameter $\sigma >0$.
  In practice however, this additional flexibility is very useful to adapt the
  model to the data, for instance by choosing on the validation set an optimal
  value for $\sigma$ with a Gaussian kernel.
\end{remark}

\section{Applications}\label{applications}

We present, in this section, several applications of the functional SVM models
described before to real world data. The first two applications illustrate the
consistent methodology introduced in section \ref{subsectionLearningSVM}: one
has an input variable with a high number of discretization points and the
second have much less discretization points. Those applications show that more
benefits are obtained from the functional approach when the data can be
reasonably considered as functions, that is when the number of discretization
points is higher than the number of observations.

The last application deals with spectrometric data and allows to show how a
functional transformation (derivative calculation) can improve the efficiency
of SVMs. For this application, we do not use the consistent methodology but a
projection on a spline space that permits easy derivative calculations.

For simplicity reasons, the parameter $C$ is chosen among a finite set of
values (in general less than 10 values) growing exponentially (for instance
0.1, 1, 10, \ldots). In each simulation, the kernel family is fixed (e.g.,
Gaussian kernels). A finite set of fully specified candidate kernels are
chosen in this family (for instance approximately 10 values of $\sigma$ in the
case of the Gaussian kernel family) and the best kernel is selected as
described in the previous section.

\subsection{Speech recognition}\label{app_biau}
We first illustrate in this section the consistent learning procedure given in
section \ref{sectionConsistency}. We compare it to the original procedure
based on $k$-nn described in \cite{BiauEtAl2005FunClassif}. In practice, the
only difference between the approaches is that we use a SVM whereas
\cite{BiauEtAl2005FunClassif} uses a $k$-nn.

The problems considered in \cite{BiauEtAl2005FunClassif} consist in
classifying speech samples\footnote{Data are available at
  \url{http://www.math.univ-montp2.fr/~biau/bbwdata.tgz}}. There are three
problems with two classes each: 
classifying ``yes'' against ``no'', ``boat'' against ``goat'' and ``sh''
against ``ao''. For each problem, we have 100 functions. Table
\ref{tableDatasetSizesBiau} gives the sizes of the classes for each problem.
\begin{table}[htbp]
  \centering
  \begin{tabular}{|c|c|c|}\hline
Problem   & Class 1 & Class $-1$ \\\hline
yes/no    & 48      & 52 \\\hline
boat/goat & 55      & 45 \\\hline
sh/ao     & 42      & 58 \\\hline
  \end{tabular}
  \caption{Sizes of the classes}
  \label{tableDatasetSizesBiau}
\end{table}

Each function is described by a vector in $\mathbb{R}^{8192}$ which
corresponds to a digitized speech frame. The goal of this benchmark is to
compare data processing methods that make minimal assumptions on the data: no
prior knowledge is used to preprocess the data.

In order to directly compare to results from \cite{BiauEtAl2005FunClassif},
performances of the algorithms are assessed by a leave-one-out 
procedure: 99 functions are used as the learning set (to which the split
sample procedure is applied to choose SVM) and the remaining
function provides a test example.

While the procedure described in \ref{subsectionLearningSVM} allows to choose
most of the parameters, both the basis $\{\Psi_j\}_{j\geq 1}$ and the penalty
term $\lambda_d$ can be freely chosen. To focus on the improvement provided by
SVM over $k$-nn, we have used the same elements as
\cite{BiauEtAl2005FunClassif}. As the data are temporal patterns,
\cite{BiauEtAl2005FunClassif} relies on the Fourier basis (moreover, the Fast
Fourier Transform allows an efficient calculation of the coordinates of the
data on the basis). The penalty term is $0$ for all $d$ below $100$ and a high
value (for instance 1000) for $d>100$. This allows to only evaluate the models
for $d\leq 100$ because the high value of $\lambda_d$ for higher $d$ prevents
the corresponding models to be chosen, regardless of their performances. As
pointed out in \cite{BiauEtAl2005FunClassif}, this choice appears to be safe
as most of the dimensions then selected are much smaller than 50.

The last free parameter is the split between the training set and the
validation set. As in \cite{BiauEtAl2005FunClassif} we have used the first 50
examples for training and the remaining 49 for validation. We report the error
rate for each problem and several methods in tables \ref{tableResultBiau} and
\ref{tableResultSVM}.
\begin{table}[htbp]
  \begin{center}
  \begin{tabular}{|c|c|c|c|c|}\hline
Problem & k-nn & QDA \\\hline
yes/no & 10\%& 7\% \\\hline
boat/goat & 21\%& 35\% \\\hline
sh/ao & 16\%& 19\%\\\hline
  \end{tabular}
\end{center}\centering
\caption{Error rate for reference methods (leave-one out)}
\label{tableResultBiau}
\end{table}

\begin{table}[htbp]
  \begin{center}
  \begin{tabular}{|c|c|c|c|}\hline
Problem/Kernel & linear (direct) & linear (projection) & Gaussian (projection)\\\hline
yes/no & 58\% & 19\%& 10\%\\\hline
boat/goat &46\%& 29\%&8\%\\\hline
sh/ao & 47\%&25\%&12\%\\\hline
  \end{tabular}
\end{center}\centering
\caption{Error rate for SVM based methods (leave-one out)}
\label{tableResultSVM}
\end{table}

Table \ref{tableResultBiau} has been reproduced from
\cite{BiauEtAl2005FunClassif}.  QDA corresponds to Quadratic Discriminant
Analysis performed, as for $k$-nn, on the projection of the data onto
a finite dimensional subspace induced by the Fourier basis. Table
\ref{tableResultSVM} gives results obtained with SVMs. The second column,
``linear (direct)'', corresponds corresponds to the direct application of 
the procedure described in \ref{subsubsectionSoftMarginSVM}, without any prior
projection. This is in fact the plain linear SVM directly applied to the
original data. The two other columns corresponds to the SVM applied to the
projected data, as described in section \ref{subsectionLearningSVM}. 

The most obvious fact is that the plain linear kernel gives very poor
performances, especially compared to the functional kernels on projections: its
results are sometimes worse than the rule that affects any observation to the
dominating class. This shows that the ridge regularization of problem
$(R_\lambda)$ is not adapted to functional data, a fact that was already known
in the context of linear discriminant analysis
\cite{HastieBujaTibshirani1995}. The projection operator improves the results
of the linear kernel, but not enough to reach the performance levels of
$k$-nn. It seems that the projected problem is therefore non linear. 

As expected, the functional Gaussian SVM performs generally better than $k$-nn
and QDA, but the training times of the methods are not comparable. On a mid
range personal computer, the full leave-one-out evaluation procedure applied
to Gaussian SVM takes approximately one and half hour (using LIBSVM
\cite{libSVM} embedded in the package e1071 of the R software
\cite{RProject}), whereas the same procedure takes only a few minutes for
$k$-nn and QDA. 

The performances of SVM with Gaussian kernel directely used on the raw data
(in \R{8192}) are not reported here as they are quite meaningless. The results
are indeed extremely sensitive to the way the grid search is conducted,
especially for the value of $C$, the regularization parameter. On the
``yes/no'' data set for instance, if the search grid for $C$ contains only
values higher than $1$, then the leave-one-out gives 19\% of error. But in each case, the
value $C=1$ is selected on the validation set. When the grid search is
extended to smaller values, the smallest value is always selected and the
error rate increases up to 46\%. Similar behaviors occur for the other data
sets. On this benchmark, the performances depend in fact on the choice of the
search grid for $C$. This is neither the case of the linear kernel on raw
data, nor the case for the projection based kernels. This is not very
surprising as Gaussian kernels have some locality problems in very high
dimensional spaces (see \cite{FrancoisASMDA2005}) that makes them difficult to
use. 

\subsection{Using wavelet basis}\label{app_hastie}
In order to investigate the limitation of the direct use of the linear SVM, we
have applied them to another speech recognition problem. We studied a part of
TIMIT database which was used in
\cite{HastieBujaTibshirani1995}\footnote{Data are available at
  \url{http://www-stat.stanford.edu/~tibs/ElemStatLearn/datasets/phoneme.data}}.
The data are log-periodograms corresponding to recording phonemes of 32 ms
duration (the length of each log-periodogram is 256). We have chosen to
restrict ourselves to classifying ``aa'' against ``ao'', because this is the
most difficult sub-problem in the database. The database is a multi-speaker
database. There are 325 speakers in the training set and 112 in the test set.
We have 519 examples for ``aa'' in the training set (759 for ``ao'') and 176
in the test set (263 for ``ao''). We use the split sample approach to choose
the parameters on the training set (50\% of the training examples are used for
validation) and we report the classification error on the test set.

Here, we do not use a Fourier basis as the functions are already represented
in a frequency form. As the data are very noisy, we decided to use a
hierarchical wavelet basis (see e.g., \cite{Mallat1989AMS}). We used the same
penalty term as in \ref{app_biau}. The error rate on the test set is reported
in table \ref{tableResultHastie}. 
\begin{table}[htbp]
  \begin{center}
  \begin{tabular}{|c|c|c|}\hline
Functional Gaussian  SVM & Functional linear SVM & Linear SVM\\\hline
22\% & 19.4\% & 20\% \\\hline
  \end{tabular}
\end{center}\centering
\caption{Error rate for all methods on the test set}
\label{tableResultHastie}
\end{table}
It appears that functional kernels are not as useful here as in the previous
example, as a linear SVM applied directly to the discretized functions (in
$\mathbb{R}^{256}$) performs as well as a linear SVM on the wavelet
coefficients. A natural explanation is that the actual dimension of the input
space (256) is smaller than the number of training examples (639) which means
that evaluating the optimal coefficients of the SVM is less difficult than in
the previous example. Therefore, the additional regularization provided by
reducing the dimension with a projection onto a small dimensional space is not
really useful in this context.

\subsection{Spectrometric data set}\label{app_spectro}

\begin{figure}[htbp]
\centering
\includegraphics[angle=270,width=0.48\textwidth]{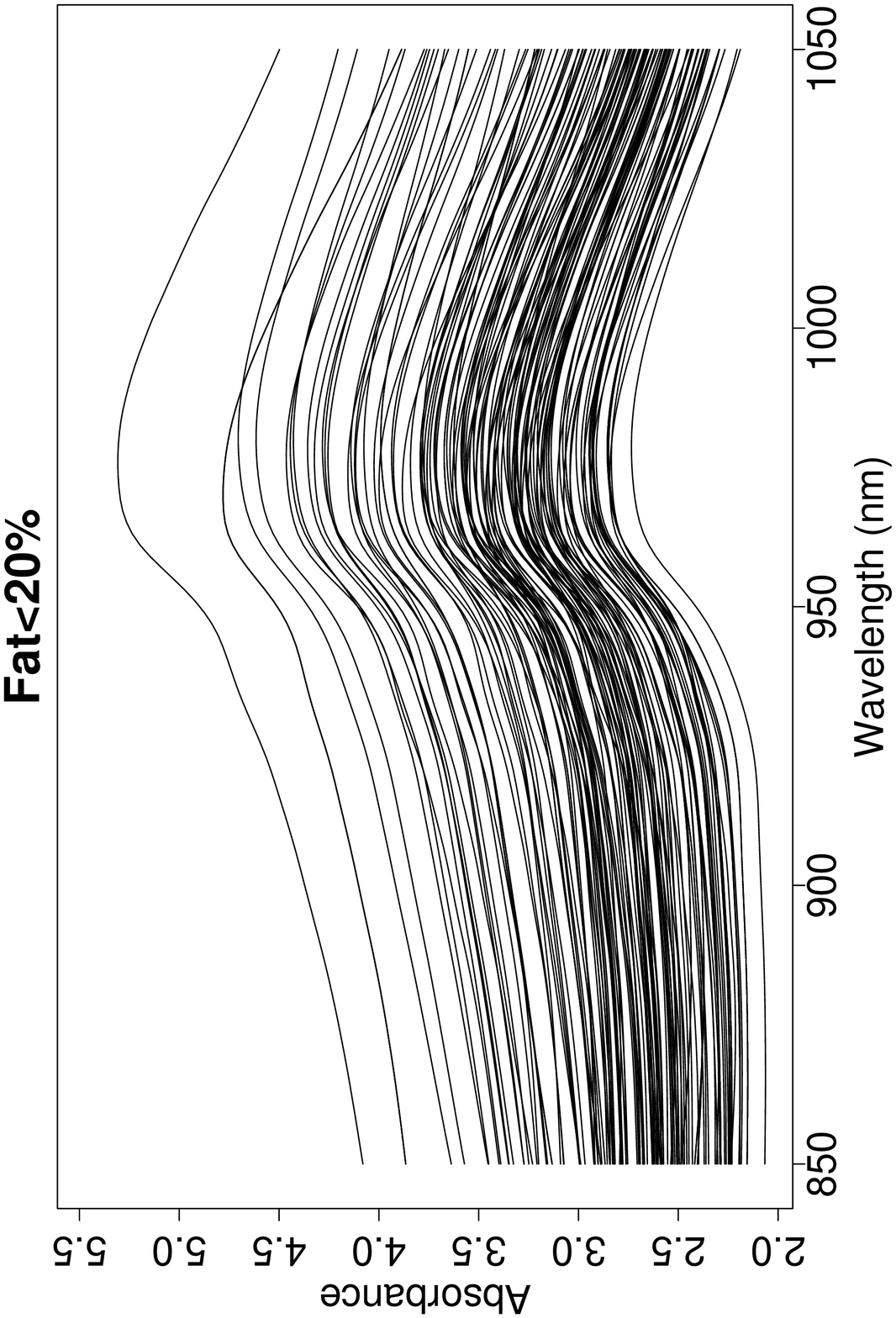}
\hfill 
\includegraphics[angle=270,width=0.48\textwidth]{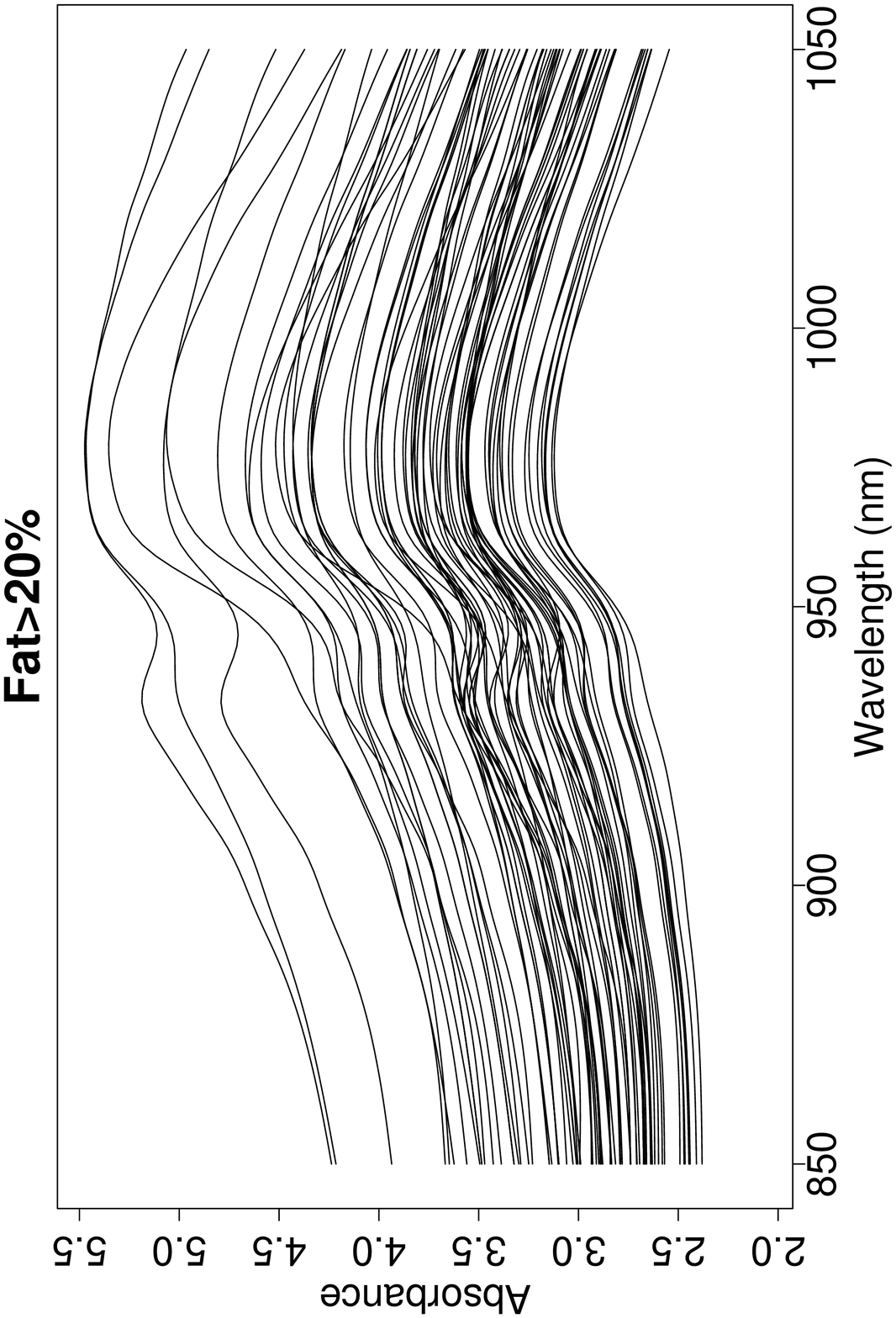}
\caption{Spectra for both classes}\label{esann_donnees spectro}
\end{figure}

We study in this section spectrometric data from food
industry\footnote{Data are available on statlib at
  \url{http://lib.stat.cmu.edu/datasets/tecator}}. Each observation is the
near infrared absorbance spectrum of a meat sample (finely chopped), recorded
on a Tecator Infratec Food and Feed Analyser (we have 215 spectra). More
precisely, an observation consists in a 100 channel spectrum of absorbances in
the wavelength range 850--1050 nm (see figure \ref{esann_donnees spectro}).
The classification problem consists in separating meat samples with a high fat
content (more than 20\%) from samples with a low fat content (less than 20\%).

\begin{figure}[htbp]
\centering
\includegraphics[angle=270,width=0.48\textwidth]{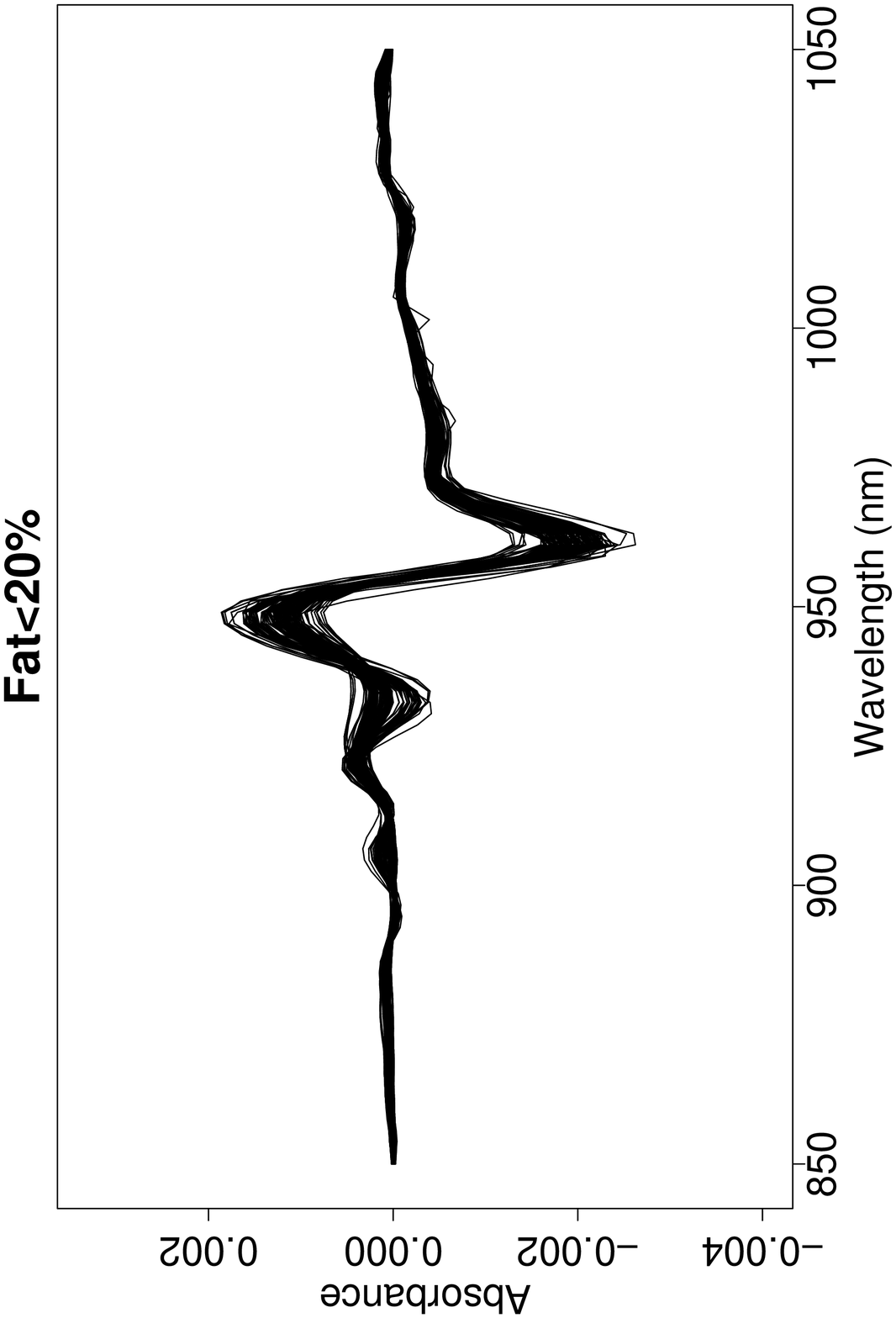}
\hfill 
\includegraphics[angle=270,width=0.48\textwidth]{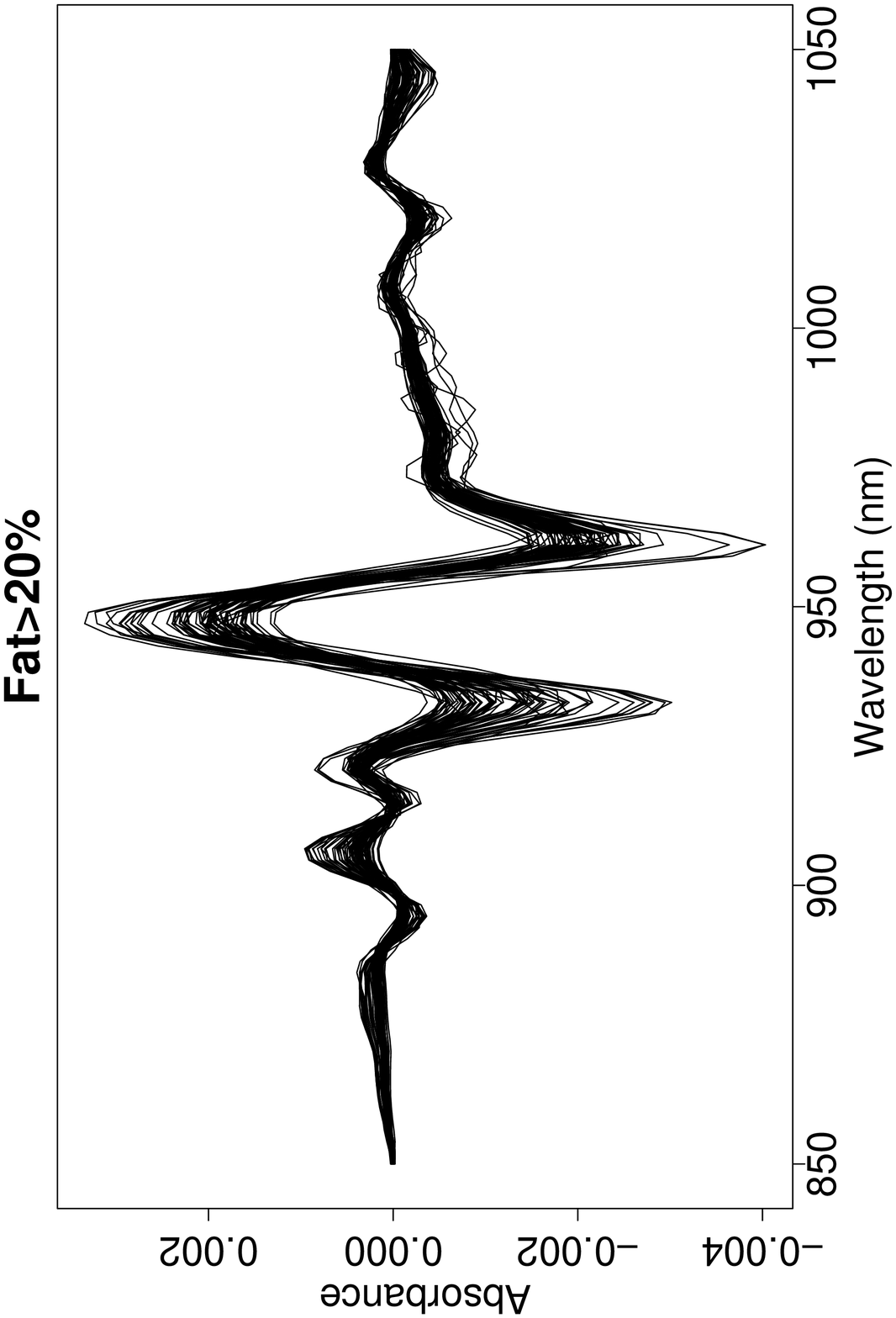}
\caption{Second derivatives of the spectra for both classes}\label{esann_donnees spectro d2}
\end{figure}

It appears on figure \ref{esann_donnees spectro} that high fat content spectra
have sometimes two local maxima rather than one: we have therefore decided to
focus on the curvature of the spectra, i.e., to use the second derivative. The
figure \ref{esann_donnees spectro d2} shows that there is more differences
between the second derivatives of each class than between the original
curves. 

The data set is split into 120 spectra for learning and 95 spectra for
testing. The problem is used to compare standard kernels (linear and Gaussian
kernels) to a derivative based kernel. We do not use here the consistent
procedure as we choose a fixed spline subspace to represent the functions so
as to calculate their second derivative. However, the parameters $C$ and
$\sigma$ are still chosen by a split sample approach that divides the 120
learning samples into 60 spectra for learning and 60 spectra for validation.
The dimension of the spline subspace is obtained thanks to a leave-one-out
procedure applied to the whole set of input functions, without taking into
account classes (see \cite{RossiEtAl05Neurocomputing} for details).

The performances depend of course on the random split between learning and
test. We have therefore repeated this splitting 250 times (as we do not select
an optimal projection dimension, the procedure is much faster than the one
used for both previous experiments). Table \ref{tableSpectro} gives the mean
error rate of those experiments on the test set.

\begin{table}[htbp]
  \begin{center}
  \begin{tabular}{|l|c|}\hline
Kernel & mean test error\\\hline\hline
Linear & 3.38\%\\\hline
Linear on second derivatives & 3.28\%\\\hline
Gaussian &7.5\%\\\hline
Gaussian on second derivatives &2.6\%\\\hline
      \end{tabular}
    \end{center}\centering
    \caption{Mean test error rate for all methods}\label{tableSpectro}
\end{table}
The results show that the problem is less difficult that the previous ones.
Nevertheless, it also appears that a functional transformation improves the
results: the use of a Gaussian kernel on second derivatives gives
significantly better results than the use of an usual kernel (linear or
Gaussian) on the original data ($t$-test results). The relatively bad
performances of the Gaussian kernel on plain data can be explained by the fact
that a direct comparison of spectra based on their $L^2(\mu)$ norm is in
general dominated by the mean value of those spectra which is not a good
feature for classification in spectrometric problems. The linear kernel is
less sensitive to this problem and is not really improved by the derivative
operator. In the Gaussian case, the use of a functional transformation
introduces expert knowledge (i.e., curvature is a good feature for some
spectrometric problems) and allows to overcome most of the limitations of the
original kernel.

\section{Conclusion}
In this paper, we have shown how to use Support Vector Machines (SVMs) for
functional data classification. While plain linear SVMs could be used directly
on functional data, we have shown the benefits of using adapted functional
kernels. We have indeed define projection based kernels that provide a
consistent learning procedure for functional SVMs. We have also introduced
transformation based kernels that allow to take into account expert knowledge
(such as the fact that the curvature of a function can be more discriminant
than its values in some applications). Both type of kernels have been tested
on real world problems. The experiments gave very satisfactory results and
showed that for some types of functional data, the performances of SVM based
classification can be improved by using kernels that make use of the
functional nature of the data.

\section*{Acknowledgements}
The authors thank the anonymous referees for their valuable suggestions that
help improving this paper.

\appendix
\section{Proofs}\label{preuves}
In order to simplify the notations, we denote $l=l_N$ when $N$ is obvious. We
also denote $X^{(d)}=P_{V_d}(X)$ and $x^{(d)}_i=P_{V_d}(x_i)$. 

The proof of the consistency result of \cite{BiauEtAl2005FunClassif} is based
on an oracle. We demonstrate a similar inequality: for $N$ large enough,
\begin{multline}\label{demo2_oracle}
 L\classifier_{a^*} - L^* \leq \inf_{d\geq 1}\left[L_d^*-L^* + \inf_{C\in\mathcal{I}_d,\ K\in \mathcal{J}_d} (L \classifier_{a} - L^*_d) + \frac{\lambda_d}{\sqrt{m}}\right]\\
+\sqrt{\frac{32(l+1)\log m}{m}} + 128 \Delta \sqrt{\frac{1}{32m(l+1)\log m}}
\end{multline}
where $m=N-l$, $\Delta \equiv \sum_{d\geq 1} |\mathcal{J}_d|
e^{-\lambda_d^2/32} < +\infty$ and $L_d^*$ is the Bayes error for the
projected problem, i.e. 
$L_d^*=\inf_{\classifier:\mathbb{R}^d\rightarrow\{-1;1\}}
\mathbb{P}(\classifier(X^{(d)})\neq Y)$.

Following \cite{BiauEtAl2005FunClassif}, we see that the definition of
$a^*=(d^*,K^*,C^*)$ leads to,
\[
\widehat{L}\classifier_{a^*} + \frac{\lambda_{d^*}}{\sqrt{m}} \leq \widehat{L} \classifier_a + \frac{\lambda_d}{\sqrt{m}}
\]

for all $a=(d,C,K)$ in $\mathcal{A}=\bigcup_{d\geq 1}\{d\}\times \mathcal{J}_d
\times [0,\mathcal{C}_d]$. Then, for all $\epsilon >0$,
\begin{multline}
   \mathbb{P}\left(L\classifier_{a^*} - \widehat{L}\classifier_a >\frac{\lambda_d}{\sqrt{m}}+\epsilon\right) \leq \mathbb{P} \left(L\classifier_{a^*} - \widehat{L}\classifier_{a^*} > \frac{\lambda_{d^*}}{\sqrt{m}}+\epsilon\right)\\
\leq \sum_{d\geq 1}\mathbb{P}\left(L\classifier_{(d,C^*,K^*)} - \widehat{L} \classifier_{(d,C^*,K^*)} >\frac{\lambda_d}{\sqrt{m}}+\epsilon\right)\\
  \leq \sum_{d\geq 1,\ K\in \mathcal{J}_d} \mathbb{P}\left(L\classifier_{(d,C^*,K)} - \widehat{L}\classifier_{(d,C^*,K)} > \frac{\lambda_d}{\sqrt{m}}+\epsilon\right)\label{demo2_etape1}
\end{multline}
In \cite{BiauEtAl2005FunClassif}, the right part of the inequality is bounded
by the use of the union bound on $\mathcal{A}$. Here, $[0,\mathcal{C}_d]$ is
not countable and thus we can not do the same. We will then use the
generalization capability of a set of linear classifiers via its shatter
coefficient. Actually, when $d$ and $K$ are set, $\classifier_{(d,C^*,K)}$ is
an affine discrimination function built from the observation projections and the kernel
$K$. More precisely, we have: 
\[
\textrm{for all } x \textrm{ in }
\mathcal{X},\qquad \classifier_{a} (x^{(d)}) = \sum_{n=1}^l \alpha_n^* y_n
K(x_n^{(d)},x^{(d)}) + b^*.  
\]
Then, $\classifier_a$ has the form
$b+\classifier$ where $\classifier$ is chosen in the set of functions spanned
by $\{K(x_1^{(d)},.),\ldots,K(x_l^{(d)},.)\}$. Let us denote by
$\mathcal{F}_K(x_1^{(d)},\ldots,x_l^{(d)})$ this set of classifiers and, for
all $\classifier$ in $\mathcal{F}_K(x_1^{(d)},\ldots,x_l^{(d)})$, we introduce
$L^l \classifier=\mathbb{P}(\classifier(X^{(d)})\neq
Y\mid(x_1,y_1),\ldots,(x_l,y_l))$. By Theorem~12.6 in
\cite{DevroyeEtAl1996Pattern}, we then have, for all $\nu>0$,
\begin{multline*}
  \mathbb{P}\left(\left. \sup_{\classifier\in
        \mathcal{F}_K(x_1^{(d)},\ldots,x_l^{(d)})}
      |\widehat{L}\classifier-L^l\classifier| >
      \nu\right|(x_1,y_1),\ldots,(x_l,y_l)\right)\leq\\
  8\mathcal{S}(\mathcal{F}_K(x_1^{(d)},\ldots,x_l^{(d)}),m)e^{-m\nu^2/32},
\end{multline*}
where $\mathcal{S}(\mathcal{F}_K(x_1^{(d)},\ldots,x_l^{(d)}),m)$ is the
shatter coefficient of $\mathcal{F}_K(x_1^{(d)},\ldots,x_l^{(d)})$, that is
the maximum number of different subsets of $m$ points that can be separated by
the set of classifiers $\mathcal{F}_K(x_1^{(d)},\ldots,x_l^{(d)})$. This set
is a vector space of dimension less or equal to $l+1$, therefore according to
chapter 13 of \cite{DevroyeEtAl1996Pattern},
$\mathcal{S}(\mathcal{F}_K(x_1^{(d)},\ldots,x_l^{(d)}),m) \leq m^{l+1}$.
This implies that, for all $(d,K)\in\mathbb{N}^*\times\mathcal{J}_d$,
\begin{eqnarray}
\lefteqn{\mathbb{P}\left(L\classifier_{(d,C^*,K)} - \widehat{L}\classifier_{(d,C^*,K)}
  > \frac{\lambda_d}{\sqrt{m}}+\epsilon\right)}&&\nonumber\\
&=&\mathbb{E}\left[\mathbb{P}\left(\left.L\classifier_{(d,C^*,K)} - \widehat{L}\classifier_{(d,C^*,K)} > \frac{\lambda_d}{\sqrt{m}}+\epsilon\right| (x_1,y_1),\ldots,(x_l,y_l)\right)\right]\nonumber\\
&\leq& \mathbb{E}\left[\mathbb{P}\left(\left.\sup_{\classifier \in
        \mathcal{F}_K(x_1^{(d)},\ldots,x_l^{(d)})} |\widehat{L}\classifier -
      L^l\classifier|
      >\frac{\lambda_d}{\sqrt{m}}+\epsilon\right|(x_1,y_1),\ldots,(x_l,y_l)\right)\right]\nonumber\\
&\leq& 8m^{l+1}e^{-\lambda_d^2/32}e^{-m\epsilon^2/32}.\label{demo2_etape2}
\end{eqnarray}
Combining (\ref{demo2_etape1}) and (\ref{demo2_etape2}), we finally see that
\[
\mathbb{P}\left(L\classifier_{a^*} - \widehat{L}\classifier_a >\frac{\lambda_d}{\sqrt{m}}+\epsilon\right)\leq 8\Delta m^{l+1}e^{-m\epsilon^2/32}.
\]
If $Z$ is a positive random variable, we have obviously
\[
\mathbb{E}(Z)\leq \mathbb{E}(Z \indic{Z>0}) = \int_0^{+\infty}\mathbb{P}(Z\geq \epsilon)\,d\epsilon.
\]
For
$Z=L\classifier_{a^*}-\widehat{L}\classifier_a-\frac{\lambda_d}{\sqrt{m}}$,
this leads, for all $a$ in
$\cup_d\{d\}\times\mathcal{I}_d\times\mathcal{J}_d$, to 
\[
L\classifier_{a^*} \leq \mathbb{E} (\widehat{L}\classifier_a) + \frac{\lambda_d}{\sqrt{m}} +\int_0^{+\infty} \mathbb{P}\left(L\classifier_{a^*} - \widehat{L}\classifier_a >\frac{\lambda_d}{\sqrt{m}}+\epsilon\right)\, d\epsilon.
\]
Finally, following \cite{BiauEtAl2005FunClassif}, for all $u>0$,
\begin{multline*}
\int_0^{+\infty} \mathbb{P}\left(L\classifier_{a^*} - \widehat{L}\classifier_a >\frac{\lambda_d}{\sqrt{m}}+\epsilon\right)\,d\epsilon \leq \int_0^u 1\, d\epsilon + \int_u^{+\infty} 8\Delta m^{l+1}e^{-m\epsilon^2/32}\,d\epsilon \\
\leq u + 128 \Delta m^{l+1} \int_u^{+\infty} \left(\frac{1}{16}+\frac{1}{m\epsilon^2}\right) e^{-m\epsilon^2/32}\, d\epsilon
\end{multline*}
and then
\[
L\classifier_{a^*} \leq \mathbb{E} (\widehat{L}\classifier_a) + \frac{\lambda_d}{\sqrt{m}} + u + \frac{128 \Delta m^l}{u}e^{-mu^2/32};
\]
if we set $u = \sqrt{\frac{32(l+1)\log m}{m}}$ and by the equality $\mathbb{E}
(\widehat{L}\classifier_a) =L\classifier_a$, we deduce that, for all $a$ in
$\mathcal{A}$, 
\[
L\classifier_{a^*} \leq L\classifier_a + \frac{\lambda_d}{\sqrt{m}} + \sqrt{\frac{32(l+1)\log m}{m}}+128\Delta\sqrt{\frac{1}{32(l+1)\log m}}
\]
which finally proves oracle (\ref{demo2_oracle}).

We conclude thanks to the following steps:
\begin{enumerate}
\item $\lim_{m\rightarrow+\infty} \sqrt{\frac{32(l+1)\log
      m}{m}}+128\mathrm{\Delta} \sqrt{\frac{1}{32m(l+1)\log m}}=0$ from the
  assumptions of Theorem~\ref{th consist};
\item Lemma 5 in \cite{BiauEtAl2005FunClassif} shows that
  $L^*_d-L^*\xrightarrow{d\rightarrow+\infty} 0$;
\item Let $\epsilon >0$. If we take a $d_0$ such that, for all $d\geq d_0$,
  $L^*_d-L^* \leq \epsilon$. To conclude, we finally have to prove that
  \[\inf_{(C,K) \in \mathcal{I}_{d_0}\times\mathcal{J}_{d_0}}
  L\classifier_{(d_0,C,K)}-L^*_{d_0} \xrightarrow{N\rightarrow +\infty} 0.\]
  This is a direct consequence of Theorem~2 in \cite{SteinwartJC2002}. Let us
  show that the hypotheses of this theorem are fulfilled:

\begin{enumerate}
  \item Theorem~2 in \cite{SteinwartJC2002} is valid for universal kernels
    that satisfy some requirements on their covering numbers. 

    As we focus on $\inf_{(C,K) \in \mathcal{I}_{d_0}\times\mathcal{J}_{d_0}}
    L\classifier_{(d_0,C,K)}$, we can choose freely the kernel and the
    regularization parameter in
    $\mathcal{I}_{d_0}\times\mathcal{J}_{d_0}$. Therefore, we choose
    $K_{d_0}$ an universal kernel with covering number of the form
    $\mathcal{O}(\epsilon^{-\nu_{d_0}})$ for some $\nu_{d_0}>0$ (this is
    possible according to our hypotheses).

  \item Theorem~2 in \cite{SteinwartJC2002} asks for $X^{(d)}$ to take its
    values in a compact set of $\mathbb{R}^d$.

  Actually, $X$ is bounded in $\mathcal{X}$ so, by definition of
  $x\rightarrow x^{(d)}$, $X^{(d)}$ takes its values in a bounded set of
  $\mathbb{R}^d$ which is included in a compact set of $\mathbb{R}^d$;

\item Finally, Theorem~2 in \cite{SteinwartJC2002} requests a particular
  behavior for $C_l$, the regularization parameter used for $l$ examples:
  $C_l$ is such that $lC_l \rightarrow +\infty$ and
  $C_l=\mathcal{O}(l^{\beta-1})$ for some $0<\beta<\frac{1}{\nu_{d_0}}$.

  Let $\beta_{d_0}$ be any number in $\left]0,\frac{1}{\nu_{d_0}} \wedge
    1\right[$ (where $a\wedge b$ denotes the infimum between $a$ and $b$).
  Then, let $C_l$ be $l^{\beta_{d_0}-1}$. This defines a sequence of real
  numbers included in $]0,1[$ which fulfills the requirements stated above. As
  $\mathcal{C}_{d_0}\geq 1$ for all $l\geq 2$, we have $C_l\in
  [0,\mathcal{C}_{d_0}]$ therefore such choice of the regularization
  parameters is compatible with the hypothesis of our theorem.
\end{enumerate}
This allows to apply Theorem~2 in \cite{SteinwartJC2002} which implies that
$L\classifier_{(d_0,(C_l),K_{d_0})}$ converges to $L^*_{d_0}$ and finally to
obtain the conclusion. 
\end{enumerate}

\bibliographystyle{abbrv}
\bibliography{total}

\end{document}